\documentclass[10pt,a4paper]{article}
\usepackage[a4paper]{geometry}
\usepackage{amssymb,latexsym,amsmath,amsfonts,amsthm}
\usepackage{graphicx}

\usepackage{epsfig}

\newtheorem{Theorem}{Theorem}[section]
\newtheorem{Definition}[Theorem]{Definition}

\newtheorem{remark}[Theorem]{Remark}
\newtheorem{Corollary}[Theorem]{Corollary}
\newtheorem{Lemma}[Theorem]{Lemma}

\newenvironment{pf}{\noindent{\bf
Proof}.\enspace}{\hfill  $\Box$ \medskip}

\numberwithin{equation}{section}

\title{Direct and Inverse Results for Multipoint Hermite-Pad\'{e} Approximants}

\date{\today}

\author{N. Bosuwan\thanks{The research of N. Bosuwan was supported by the Strengthen Research Grant for New Lecturer from the Thailand Research Fund and the Office of the Higher Education Commission (MRG6080133) and Faculty of Science, Mahidol University.}, G. L\'{o}pez Lagomasino\footnotemark[2], Y. Zaldivar Gerpe\footnotemark[2]}

\begin{document}

\maketitle
\renewcommand{\thefootnote}{\fnsymbol{footnote}}
\footnotetext[2]{
G. L\'{o}pez Lagomasino and Y. Zaldivar Gerpe received support from research grant MTM 2015-65888-C4-2-P of Ministerio de Econom\'{\i}a, Industria y Competitividad, Spain.}

\begin{abstract}
Given a system of functions $\mathbf{f}=(f_1,\ldots,f_d)$ analytic on a neighborhood of some compact subset $E$ of the complex plane, we give necessary and sufficient conditions for the convergence with geometric rate of the common denominators of multipoint Hermite-Pad\'{e} approximants. The exact rate of convergence of these denominators and of the approximants themselves is given in terms of the analytic properties of the system of functions. These  results allow to detect the location of the poles of the system of functions which are in some sense ``closest'' to $E$.

\end{abstract}

\medskip


\medskip


\maketitle

\section{Statement of the main result.}

We shall consider a general interpolation scheme for constructing vector rational approximations to a given vector of analytic functions which generalizes the construction of the classical Hermite-Pad\'{e} approximants.

Let $E$ be a bounded continuum with connected complement in the complex plane $\mathbb{C}$. By $\mathcal{H}(E)$ we denote the space of all functions holomorphic in some neighborhood of $E$. Set
\[
\mathcal{H}(E)^d:=\lbrace(f_1,\dots,f_d):f_j\in \mathcal{H}(E),j=1,\dots,d\rbrace.
\]
Let $\alpha \subset E$ be a table of points; more precisely, $\alpha=\lbrace \alpha_{n,k}\rbrace$, $k=1,\ldots,n$, $n=1,2,\ldots$.  We propose the following definition.

\begin{Definition}\label{multipdef} \textup{
Let ${\bf f}\in \mathcal{H}(E)^d$. Fix a multi-index ${\bf m}=(m_1,\dots,m_d)\in\mathbb{N}^d$ and $n\in\mathbb{N}$. Set $|{\bf m}|=m_1+\cdots+m_d$. Then, there exist polynomials  $Q_{n,{\bf m}}$, $P_{n,{\bf m},k}$, $k=1,\ldots,d$ such that
\begin{enumerate}
\item[b.1)] $\deg P_{n,{\bf m},k}\leq n-m_k$, $\deg Q_{n,{\bf m}}\leq |\mathbf{m}|$, $Q_{n,{\bf m}}\not\equiv 0$,
\item[b.2)] ${(Q_{n,{\bf m}}f_k-P_{n,{\bf m},k})}/{a_{n+1}}\in \mathcal{H}(E),$
\end{enumerate}
where $a_n(z)=\prod_{k=1}^{n}(z-\alpha_{n,k})$. The vector rational function
\[
\mathbf{R}_{n,\mathbf{m}}=\left(R_{n,{\bf m},1},\dots,R_{n,{\bf m},d}\right)={\left( P_{n,{\bf m},1},\dots,P_{n,{\bf m},d}\right)}/{Q_{n,{\bf m}}}
\]
is called a \emph{multipoint Hermite-Pad\'e (MHP) approximant of $\mathbf{f}$ with respect to ${\bf m}$ and $\alpha$}.
}
\end{Definition}

This vector rational approximation, in general, is not uniquely determined. Hereafter, we assume that given $(n,{\bf m})$, one particular solution is taken. Without loss of generality we can assume that $Q_{n,{\bf m}}$ is a monic polynomial that has no common zero simultaneously with all $P_{n,{\bf m},k}$. In all what follows $\bf m$ remains fixed and $\{\mathbf{R}_{n,\mathbf{m}}\}_{n \in \mathbb{N}}$ is called a row sequence of MHP of $\bf f$ with respect to $\bf m$.

\medskip

Multipoint Hermite-Pad\'e approximation reduces to classical Hermite-Pad\'{e} approximation when $E$ is a disk about the origin and $a_n(z)=z^n$. There are not many papers dealing with the convergence properties of row sequences of Hermite-Pad\'{e} approximation. The first significant contribution in this direction is due to Graves-Morris and Saff in \cite{graves}, where an analogue of the Montessus de Ballore theorem \cite{Mon} was proved. In that paper, the authors studied the classical case and stated a result for multipoint interpolation. They assume that the system of approximated functions is, so called, polewise independent. More recently, the authors of \cite{cacoq1} and \cite{cacoq2} managed to weaken the assumption of polewise independence obtaining  sharp estimates of the rate of convergence, improving the region of convergence, and giving an analogue of Gonchar's converse statement to the Montessus de Ballore theorem  for row sequences of Pad\'e approximants (see Remark in \cite{gon1}, also \cite{gon3} and \cite{gon2}). Here, we generalize the results in \cite{cacoq2} to MHP approximants. Extensions  in other directions using expansions in  orthogonal and Faber polynomials of the vector function to produce the vector rational approximants of $\bf f$ were provided in \cite{Bos1,Bos2}. For other approaches to the study of row sequences of vector rational approximation see \cite{Sidi} and \cite{VanBarel}.

\medskip

In the study of the convergence of general interpolation schemes, it is common to impose on the table of interpolation nodes various restrictions which determine the asymptotic behavior of the sequence of polynomials $a_n$. Let $\Phi_E$ be a holomorphic univalent function mapping the complement of $E$ onto the exterior of the closed unit disk with $\Phi_E(\infty)=\infty$ and $\Phi'_E(\infty)> 0$. It is well known that there exist tables of points $\alpha$ satisfying the condition
\begin{equation}\label{raiz}
\lim_{n \to \infty} |a_n(z)|^{1/n} = c|\Phi_E(z)|,
\end{equation}
or the stronger condition
\begin{equation}\label{eq5}
\lim_{n\to \infty}a_n(z)/c^n\Phi_E^n(z) = G(z)\neq 0,
\end{equation}
uniformly on compact subsets of $\overline{\mathbb{C}}\setminus E$, where $c$ denotes some positive constant, see \cite[Chapters 8-9]{Walsh}. For each $\rho>1$, we introduce
\[
\Gamma_{\rho}:=\lbrace z\in\mathbb{C}: |\Phi_E(z)|=\rho\rbrace,\quad\text{and}\quad D_{\rho}:=E\cup\lbrace z\in\mathbb{C}: |\Phi_E(z)|<\rho\rbrace
\]
as the \emph{level curve of index $\rho$} and the \emph{canonical domain of index $\rho$}, respectively. Let $\rho_0({\bf f})$ be equal to the index $\rho$ of the largest canonical domain $D_\rho$ to which all $f_k$, $k=1,\dots,d$ can be extended as holomorphic functions simultaneously.

\medskip

Gonchar proved the following analogue of the Cauchy-Hadamard formula for $f \in \mathcal{H}(E)$ and interpolation tables satisfying \eqref{eq5}:
\begin{equation}\label{cauchyhadamard}
\rho_0(f)= \left(c \cdot  \limsup_{n\rightarrow\infty}\left|\int_{\Gamma}\frac{f(t)}{a_{n+1}(t)}dt\right|^{1/n} \right)^{-1},
\end{equation}
where $\Gamma$ is a contour encircling $E$ and lying in the domain of holomorphy of $f.$
This formula is a special case of \cite[Corollary 3]{Bus1}. (We point out that \eqref{cauchyhadamard} is displayed as formula (17) in \cite{Bus1}, but with the typo that $c$ is missing.)

\begin{Definition}\label{systempole2}\textup{
Given $\mathbf{f}=(f_1,\ldots,f_d)\in \mathcal{H}(E)^d$ and $\mathbf{m}=(m_1,\ldots,m_d)\in \mathbb{N}^d$ we say that \emph{$\xi\in\mathbb{C}$ is a system pole of order $\tau$ of $(\mathbf{f}, \mathbf{m})$} if $\tau$ is the largest positive integer such that for each $s=1,\ldots,\tau$ there exists at least one polynomial combination of the form
\begin{equation}\label{eqsystempole2}
\sum\limits_{k=1}^dp_kf_k,\qquad \deg p_k<m_k,\qquad k=1,\ldots,d,
\end{equation}
which is analytic in a neighborhood of $\overline{D}_{|\Phi_E(\xi)|}$ except for a pole at $z=\xi$ of exact order $s$.}
\end{Definition}

\medskip

The concept of system pole depends not only on the system of functions ${\bf f}$ but also on the multi index ${\bf m}$. For example, poles of the individual functions $f_k$ need not be system poles of $({\bf f,m})$ and system poles need not be poles of any of the functions $f_k$ (see examples in \cite{cacoq2}). It is easy to see that system poles also depend on $\alpha$, or more precisely on the geometry of the associated canonical regions. However, since $\bf m$ and $\alpha$ will remain fixed, occasionally we may simply refer to system poles of $\bf f$.

\medskip

Let $\tau$ be the order of $\xi$ as a system pole of ${\bf f}$. For each $s=1,\dots,\tau$, let $\rho_{\xi,s}({\bf f},{\bf m})$ denote the largest of all the numbers $\rho_s(g)$ (the index of the largest canonical domain containing at most $s$ poles of $g$), where $g$ is a polynomial combination of type \eqref{eqsystempole2} that is holomorphic on a neighborhood of $\overline{D}_{|\Phi_E(\xi)|}$ except for a pole at $z=\xi$ of order $s$. Then, we define
\[
R_{\xi,s}({\bf f},{\bf m}):=\min_{k=1,\dots,s} \rho_{\xi,k}({\bf f},{\bf m}),
\]
and
\[
R_{\xi}({\bf f},{\bf m}):=R_{\xi,\tau}({\bf f},{\bf m})=\min_{k=1,\dots,\tau} \rho_{\xi,k}({\bf f},{\bf m}).
\]

Fix $k=\lbrace 1,\dots,d\rbrace$. Let $D_{k}({\bf f},{\bf m})$ be the largest canonical domain in which all the poles of $f_k$ are system poles of ${\bf f}$ with respect to ${\bf m}$, their order as poles of $f_k$ does not exceed their order as system poles, and $f_k$ has no other singularity. By $R_{k}({\bf f},{\bf m})$, we denote the index of this canonical domain. Let $\xi_1,\dots,\xi_N$ be the poles of $f_k$ in $D_{k}({\bf f},{\bf m})$. For each $j=1,\dots,N$, let $\hat{\tau}_j$ be the order of $\xi_j$ as pole of $f_k$ and $\tau_j$ be its order as a system pole. By assumption, $\hat{\tau}_j\leq\tau_j$. Set
\[
R_{k}^*({\bf f},{\bf m}):=\min\left\lbrace R_{k}({\bf f},{\bf m}),\min_{j=1,\dots,N}R_{\xi_j,\hat{\tau}_j}({\bf f},{\bf m})\right\rbrace
\]
and let $D_{k}^*({\bf f},{\bf m})$ be the canonical domain with this index.

\medskip

By $Q_{\mathbf{m}}^{\mathbf{f}}$ we denote the monic polynomial whose zeros are the system poles of ${\bf f}$ with respect to ${\bf m}$ taking account of their order. The set of distinct zeros of $Q_{\mathbf{m}}^{\mathbf{f}}$ is denoted by $\mathcal{P}_{\mathbf{m}}^{\mathbf{f}}$.

\medskip

The following theorem constitutes our main result.

\begin{Theorem}\label{maintheorem}
Suppose \eqref{eq5} takes place. Let $\mathbf{f}\in \mathcal{H}(E)^d$ and fix a multi-index $\mathbf{m}\in\mathbb{N}^d$. Then, the next two assertions are equivalent:
\begin{enumerate}
\item[(a)] $\mathbf{f}$ has exactly $|\mathbf{m}|$ system poles with respect to $\mathbf{m}$ counting multiplicities.
\item[(b)] For all sufficiently large n, the denominators $Q_{n,\mathbf{m}}$ of multipoint Hermite-Pad\'e approximants of $\mathbf{f}$ are uniquely determined and there exists a polynomial $Q_{\mathbf{m}}$ of degree $|\mathbf{m}|$ such that
\begin{equation}\label{eq1teo}
\limsup\limits_{n\rightarrow\infty}\Vert Q_{n,\mathbf{m}} - Q_{\mathbf{m}}\Vert^{1/n}=\theta<1,
\end{equation}
\end{enumerate}
where $\|\cdot\|$ denotes the coefficient norm in the space of polynomials of degree $\leq |\bf m|$.
Moreover, if either (a) or (b) takes place, then $Q_{\mathbf{m}}\equiv Q_{\mathbf{m}}^{\mathbf{f}}$,
\begin{equation}\label{eq6}
\theta=\max \left\lbrace\frac{|\Phi_E(\xi)|}{R_{\xi}(\mathbf{f},\mathbf{m})}:\xi\in\mathcal{P}_{\mathbf{m}}^{\mathbf{f}}\right\rbrace,
\end{equation}
and for any compact subset $\mathcal{K}$ of $D_{k}^*({\bf f},{\bf m})\setminus\mathcal{P}_{\mathbf{m}}^{\mathbf{f}}$,
\begin{equation}\label{eq6.1}
\limsup_{n\rightarrow\infty}\|R_{n,{\bf m},k}-f_k\|_{\mathcal{K}}^{1/n}\leq\frac{\|\Phi_E\|_\mathcal{K}}{R_{k}^*({\bf f},{\bf m})},
\end{equation}
where $\|\cdot\|_{\mathcal{K}}$ denotes the sup-norm on $\mathcal{K}$ and if $\mathcal{K}\subset E$, then $\|\Phi_E\|_{\mathcal{K}}$ is replaced by $1$.
\end{Theorem}

\section{Direct statements}

\subsection{An auxiliary result}
For each $n\geq |\bf{m}|$, let $q_{n,{\bf m}}$ be the polynomial $Q_{n,{\bf m}}$ normalized so that
\begin{equation}\label{norm}
\sum\limits_{k=0}^{|\bf{m}|}|\lambda_{n,k}|=1,\qquad q_{n,{\bf m}}(z)=\sum\limits_{k=0}^{|\bf{m}|}\lambda_{n,k}z^k.
\end{equation}
This normalization implies that the polynomials $q_{n,{\bf m}}$ are uniformly bounded on each compact subset of $\mathbb{C}$.

\begin{Lemma} Let ${\bf f}\in \mathcal{H}(E)^d$ and fix a multi-index ${\bf m}\in \mathbb{N}^d.$ Assume that \eqref{eq5} takes place and $\xi$ is a system pole of order $\tau$ of ${\bf f}$ with respect to ${\bf m}$. Then
\begin{equation}\label{eq induc}
\limsup_{n\rightarrow\infty}|q_{n,\bf{m}}^{(s)}(\xi)|^{1/n}\leq\frac{|\Phi_E(\xi)|}{R_{\xi,s+1}(\bf{f},\bf{m})}, \qquad s=0,\ldots,\tau -1.
\end{equation}
\end{Lemma}

\begin{pf} Consider a polynomial combination $g_1$ of type \eqref{eqsystempole2} that is analytic on a neighborhood of $\overline{D}_{|\Phi_E(\xi)|}$ except for a simple pole $z=\xi$ and verifies that $\rho_1(g_1)=R_{\xi,1}({\bf f},{\bf m})(=\rho_{\xi,1}(\bf{f},\bf{m}))$. Then, we have
\[
g_1=\sum\limits_{k=1}^{d}p_{k,1}f_k,\quad\deg p_{k,1}<m_k,\quad k=1,\ldots,d.
\]
Define $h_1(z)=(z-\xi)g_1(z)$. The function
\[
\frac{q_{n,{\bf m}}(z)h_1(z)}{a_{n+1}(z)}-\frac{z-\xi}{a_{n+1}(z)}\sum\limits_{k=1}^{d}p_{k,1}(z)P_{n,{\bf m},k}(z)
\]
is analytic on $D_{\rho_1(g_1)}$. Take $1<\rho<\rho_1(g_1)$, and set $\Gamma_\rho=\lbrace z\in \mathbb{C}: |\Phi_E(z)|=\rho\rbrace$. Set $P_{n,1}(z) =  \sum\limits_{k=1}^{d}p_{k,1}(z)P_{n,{\bf m},k}(z)$. Since $\deg  (z-\xi)P_{n,1}(z) \leq n$, we have
\[ \frac{1}{2\pi i}\int_{\Gamma_\rho} \frac{(t-\xi)P_{n,1}(t)}{(t-z)a_{n+1}(t)}dt = 0.
\]
Using Hermite's interpolation formula (see \cite{Walsh}), we obtain
\[
q_{n,{\bf m}}(z)h_1(z)-(z-\xi)\sum\limits_{k=1}^{d}p_{k,1}P_{n,{\bf m},k}(z)=\frac{1}{2\pi i}\int_{\Gamma_\rho}\frac{a_{n+1}(z)}{a_{n+1}(t)}\frac{q_{n,{\bf m}}(t)h_1(t)}{t-z}dt,
\]
for all $z$ with $|\Phi_E(z)|<\rho$. In particular, taking $z=\xi$ in the above formula, we arrive at
\begin{equation}\label{eq7}
q_{n,{\bf m}}(\xi)h_1(\xi)=\frac{1}{2\pi i}\int_{\Gamma_\rho}\frac{a_{n+1}(\xi)}{a_{n+1}(t)}\frac{q_{n,{\bf m}}(t)h_1(t)}{t-\xi}dt.
\end{equation}
Then, taking account of \eqref{eq5}, it easily follows that
\[
\limsup_{n\rightarrow\infty}|q_{n,{\bf m}}(\xi)h_1(\xi)|^{1/n}\leq\frac{|\Phi_E(\xi)|}{\rho}.
\]
Using that $h_1(\xi)\neq 0$ and making $\rho$ tend to $\rho_1(g_1)$, we obtain
\[
\limsup_{n\rightarrow\infty}|q_{n,{\bf m}}(\xi)|^{1/n}\leq\frac{|\Phi_E(\xi)|}{R_{\xi,1}({\bf f},{\bf m})}<1.
\]

Now, we employ induction. Suppose that
\begin{equation}\label{eq8}
\limsup_{n\rightarrow\infty}|q_{n,{\bf m}}^{(j)}(\xi)|^{1/n}\leq\frac{|\Phi_E(\xi)|}{R_{\xi,j+1}({\bf f},{\bf m})},\qquad j=0,1,\ldots,s-2,
\end{equation}
where $s \leq \tau$. Let us prove that formula \eqref{eq8} holds for $j=s-1$. This will imply \eqref{eq induc}.

\medskip

Consider a polynomial combination $g_s$ of type \eqref{eqsystempole2} that is analytic on a neighborhood of $\overline{D}_{|\Phi_E(\xi)|}$ except for a pole of order $s$ at $z=\xi$ and verifies that $\rho_s(g_s)=R_{\xi,s}({\bf f},{\bf m})$. Then,
\[
g_s=\sum\limits_{k=1}^{d}p_{k,s}f_k,\ \ \ \deg p_{k,s}<m_k,\ \ \ k=1,\ldots,d.
\]
Set $h_s(z)=(z-\xi)^sg_s(z)$. The function
\[
\frac{q_{n,{\bf m}}(z)h_s(z)}{a_{n+1}(z)(z-\xi)^{s-1}}-\frac{z-\xi}{a_{n+1}(z)}\sum\limits_{k=1}^{d}p_{k,s}(z)P_{n,{\bf m},k}(z)
\]
is analytic on $D_{\rho_s(g_s)}\setminus\lbrace\xi\rbrace$. Set $P_{n,s}=\sum\limits_{k=1}^{d}p_{k,s}P_{n,{\bf m},k}$. Fix an arbitrary compact set $\mathcal{K}\subset D_{\rho_s(g_s)} \setminus\lbrace\xi\rbrace$. Take $\delta>0$ sufficiently small so that $\{z\in \mathbb{C}: |z-\xi|\leq \delta\}\cap \mathcal{K}=\emptyset$ and $1<\rho <\rho_s(g_s)$. Using Hermite's interpolation formula, for all $z\in\mathcal{K}$,  we have
\begin{equation}\label{eq9}
\frac{q_{n,{\bf m}}(z)h_s(z)}{(z-\xi)^{s-1}}-(z-\xi)P_{n,s}(z)=I_n(z)-J_n(z),
\end{equation}
where
\[
I_n(z)=\frac{1}{2\pi i}\int_{\Gamma_{\rho}}\frac{a_{n+1}(z)}{a_{n+1}(t)}\frac{q_{n,{\bf m}}(t)h_s(t)}{(t-\xi)^{s-1}(t-z)}dt
\]
and
\[
J_n(z)=\frac{1}{2\pi i}\int_{|t-\xi|=\delta}\frac{a_{n+1}(z)}{a_{n+1}(t)}\frac{q_{n,{\bf m}}(t)h_s(t)}{(t-\xi)^{s-1}(t-z)}dt.
\]
The first integral $I_n$ is estimated as in \eqref{eq7} to obtain
\begin{equation}\label{eq10}
\limsup_{n\rightarrow\infty}\|I_n\|_{\mathcal{K}}^{1/n}\leq\frac{\|\Phi_E\|_{\mathcal{K}}}{\rho_s(g_s)}.
\end{equation}
For $J_n$, as $\deg q_{n,{\bf m}}\leq |{\bf m}|$ write
\[
q_{n,{\bf m}}(t)=\sum\limits_{j=0}^{|{\bf m}|}\frac{q_{n,{\bf m}}^{(j)}(\xi)}{j!}(t-\xi)^j.
\]
Then
\begin{equation}\label{eq11}
J_n(z)=\sum\limits_{j=0}^{s-2}\frac{1}{2\pi i}\int_{|t-\xi|=\delta}\frac{a_{n+1}(z)}{a_{n+1}(t)}\frac{h_s(t)}{(t-\xi)^{s-1-j}}\frac{q_{n,{\bf m}}^{(j)}(\xi)}{j!(t-z)}dt.
\end{equation}
Using the induction hypothesis \eqref{eq8}, from \eqref{eq11} it easily follows that
\begin{equation}\label{eq11.1}
\limsup_{n\rightarrow\infty}\|J_n\|_{\mathcal{K}}^{1/n}\leq\frac{\|\Phi_E\|_{\mathcal{K}}}{|\Phi_E(\xi)|}\frac{|\Phi_E(\xi)|}{R_{\xi,s-1}({\bf f},{\bf m})}=\frac{\|\Phi_E\|_{\mathcal{K}}}{R_{\xi,s-1}({\bf f},{\bf m})}.
\end{equation}
Now, \eqref{eq9}, \eqref{eq10}, and \eqref{eq11.1} give
\begin{equation}\label{eq12}
\limsup_{n\rightarrow\infty}\|q_{n,{\bf m}}h_s-(z-\xi)^s P_{n,s}\|_{\mathcal{K}}^{1/n}\leq\frac{\|\Phi_E\|_{\mathcal{K}}}{R_{\xi,s}({\bf f},{\bf m})}.
\end{equation}

As the function inside the norm in \eqref{eq12} is analytic in $D_{\rho_l(g_l)}$, from the maximum principle it follows that \eqref{eq12} also holds for any compact set $\mathcal{K}\subset D_{\rho_l(g_l)}$. Using Cauchy's integral formula, from \eqref{eq12} we also obtain that
\begin{equation}\label{eq12.1}
\limsup_{n\rightarrow\infty}\|(q_{n,{\bf m}}h_s-(z-\xi)^s P_{n,s})^{(s-1)}\|_{\mathcal{K}}^{1/n}\leq\frac{\|\Phi_E\|_{\mathcal{K}}}{R_{\xi,s}({\bf f},{\bf m})}.
\end{equation}
Taking $z=\xi$ in \eqref{eq12.1}, we have
\[
\limsup_{n\rightarrow\infty}|(q_{n,{\bf m}}h_s)^{(s-1)}(\xi)|^{1/n}\leq\frac{|\Phi_E(\xi)|}{R_{\xi,s}({\bf f},{\bf m})}.
\]
Using the Leibniz formula for higher derivatives of a product of two functions, the induction hypothesis \eqref{eq8}, and that $h_s(\xi)\neq 0$, it follows that
\[
\limsup_{n\rightarrow\infty}|q_{n,{\bf m}}^{(s-1)}(\xi)|^{1/n}\leq\frac{|\Phi_E(\xi)|}{R_{\xi,s}({\bf f},{\bf m})},
\]
This completes the induction and the proof.
\end{pf}

\subsection{Proof of $(a)\Rightarrow (b)$}
Let $\lbrace\xi_1,\dots,\xi_p\rbrace$ be the distinct system poles of ${\bf f}$ with respect to ${\bf m}$, and let $\tau_j$ be the order of $\xi_j$ as a system pole, $j=1,\dots,p$. By assumption, $\tau_1+\dots +\tau_p=|{\bf m}|$. We have proved that, for $j=1,\dots,p$ and  $s=0,1,\dots,\tau_j-1,$
\begin{equation}\label{eq14}
\limsup_{n\rightarrow\infty}|q_{n,{\bf m}}^{(s)}(\xi_j)|^{1/n}\leq\frac{|\Phi_E(\xi_j)|}{R_{\xi_j,s+1}({\bf f},{\bf m})}\leq\frac{|\Phi_E(\xi_j)|}{R_{\xi_j}({\bf f},{\bf m})},
\end{equation}
where $R_{\xi_j}({\bf f},{\bf m}):=R_{\xi_j,\tau_j}({\bf f},{\bf m})$.
Using the Hermite interpolation, it is easy to construct a basis $\lbrace \ell_{j,s}\rbrace, 1\leq j\leq p, 0\leq s\leq\tau_j-1,$ in the space of polynomials of degree at most $|{\bf m}|-1$ satisfying
\[
\ell_{j,s}^{(k)}(\xi_i)=\delta_{i,j}\delta_{k,s},\qquad 1\leq i\leq p,\qquad 0\leq k\leq\tau_i-1.
\]
Then,
\begin{equation}\label{eq15}
q_{n,{\bf m}}(z)=\sum\limits_{j=1}^{p}\sum\limits_{s=0}^{\tau_j-1}q_{n,{\bf m}}^{(s)}(\xi_j)\ell_{j,s}(z)+\lambda_{n,|{\bf m}|}Q_{{\bf m}}^{{\bf f}}.
\end{equation}
Using \eqref{eq14} and \eqref{eq15}, we have for any compact set $\mathcal{K} \subset\mathbb{C}$,
\begin{equation}\label{eq15.1}
\limsup_{n\rightarrow\infty}\|q_{n,{\bf m}}-\lambda_{n,|{{\bf m}}|}Q_{{\bf m}}^{\bf f}\|_{\mathcal{K}}^{1/n}\leq\theta,
\end{equation}
 where
\begin{equation}\label{eq16}
\theta=\max\left\lbrace\frac{|\Phi_E(\xi)|}{R_{\xi}(\bf{f},{\bf{m}})}:\ \xi\in\mathcal{P}_{{\bf m}}^{\bf f}\right\rbrace<1.
\end{equation}
Now, necessarily
\begin{equation}\label{eq17}
\liminf_{n\rightarrow\infty}|\lambda_{n,|{\bf m}|}|>0.
\end{equation}
Indeed, if there is a subsequence of indices $\Lambda\subset\mathbb{N}$ such that $\lim_{n\in\Lambda}|\lambda_{n,|\bf{m}|}|=0,$ then from \eqref{eq16}, as the polynomials $q_{n,{\bf m}}$ converge, we would have that $\lim_{n\in\Lambda} q_{n,{\bf m}}=0$ which contradicts \eqref{norm}. Since
\begin{equation}\label{eq18}
q_{n,{\bf m}}=\lambda_{n,|{\bf m}|}Q_{n,{\bf m}},
\end{equation}
from \eqref{eq15.1} and \eqref{eq17} it follows that
\begin{equation}\label{eq19}
\limsup_{n\rightarrow\infty}\|Q_{n,{\bf m}}-Q_{{\bf m}}^{\bf f}\|_{\mathcal{K}}^{1/n}\leq\theta.
\end{equation}
In finite dimensional spaces all norms are equivalent; therefore, \eqref{eq19} is also true with the coefficient norm which means that \eqref{eq1teo} is satisfied with $=$ replaced by $\leq$.

\medskip

In particular, for all sufficiently large $n$ necessarily $\deg Q_{n,{\bf m}}=|\mathbf{m}|$. The difference of any two distinct monic polynomials satisfying Definition \ref{multipdef} with the same degree produces a new solution of degree strictly less than $|\mathbf{m}|$, but we have proved that any solution must have degree $|\mathbf{m}|$ for all sufficiently large $n$. Hence, the polynomial $Q_{n,{\bf m}}$ is uniquely determined for all sufficiently large $n$.

\medskip

Now, we prove the equality in \eqref{eq1teo}. To the contrary, suppose that
\begin{equation}\label{eq20}
\limsup_{n\rightarrow\infty}\|Q_{n,{\bf m}}-Q_{{\bf m}}^{\bf f}\|^{1/n}<\max\left\lbrace
\frac{|\Phi_E(\xi)|}{R_{\xi}(\bf{f},{\bf{m}})}:\ \xi\in\mathcal{P}_{{\bf m}}^{\bf f}\right\rbrace.
\end{equation}
Let $\zeta$ be a system pole of $\bf{f}$ such that
\begin{equation}\label{d}
\frac{|\Phi_E(\zeta)|}{R_{\zeta}(\bf{f},{\bf{m}})}=\max\left\lbrace
\frac{|\Phi_E(\xi)|}{R_{\xi}(\bf{f},{\bf{m}})}:\ \xi\in\mathcal{P}_{{\bf m}}^{\bf f}\right\rbrace.
\end{equation}
Clearly, the inequality \eqref{eq20} implies that $R_{\zeta}({\bf f},{{\bf m}})<\infty$.

\medskip

Choose a polynomial combination
\begin{equation}\label{eq21}
g=\sum\limits_{k=1}^{d}p_{k}f_k,\quad \deg p_{k}<m_k,\quad k=1,\dots,d,
\end{equation}
that is holomorphic on a neighborhood of $\overline{D}_{|\Phi_E(\zeta)|}$ except for a pole of some order $l$ at $z=\zeta$ with $\rho_l(g)=R_{\zeta}({\bf f},{{\bf m}})$.
Notice that $Q_{\textup{\textbf{m}}}^{\textup{\textbf{f}}}g$ must have a singularity on the boundary of $D_{\rho_l(g)}$ which implies
\begin{equation}\label{a}
\frac{1}{R_{\zeta}({\bf f},{\bf m})}=c \cdot  \limsup_{n\rightarrow\infty}\left|\int_{\Gamma_{\rho}}\frac{Q_{{\bf m}}^{\bf f}(t)g(t)}{a_{n+1}(t)}dt\right|^{1/n}.
\end{equation}
In fact, if $Q_{\textup{\textbf{m}}}^{\textup{\textbf{f}}}g$ had no singularity on the boundary of $D_{\rho_l(g)}$, then all singularities of $g$  on the boundary of $D_{\rho_l(g)}$ would be at most poles and their order as poles of $g$ would be smaller than their order as system poles of ${\textup{\textbf{f}}}$. In this case, we could find a different polynomial combination $g_1$ of type \eqref{eq21} for which $\rho_l(g_1)>\rho_l(g)=R_{\zeta}(\textup{\textbf{f}},\textup{\textbf{m}})$ which contradicts the definition of $R_{\zeta}(\textup{\textbf{f}},\textup{\textbf{m}})$. Therefore, $Q_{\textup{\textbf{m}}}^{\textup{\textbf{f}}}g$ has a singularity on the the boundary of $D_{\rho_l(g)}$ and the equality \eqref{a} holds.

Now,
\[ \left(Q_{n,{\bf m}}(z)g(z) - \sum_{k=1}^d p_k(z)P_{n,{\bf m},k}(z)\right)/a_{n+1}(z)
\]
is holomorphic in $D_{\rho_l(g)}$ and $\deg \sum_{k=1}^d p_kP_{n,{\bf m},k} < n$; therefore, from Cauchy's integral theorem we have that
\begin{equation}\label{b} 0 = \int_{\Gamma_{\rho}} \frac{Q_{n,{\bf m}}(z)g(z) - \sum_{k=1}^d p_k(z)P_{n,{\bf m},k}(z)}{a_{n+1}(z)} dz = \int_{\Gamma_{\rho}} \frac{Q_{n,{\bf m}}(z)g(z)}{a_{n+1}(z)} dz,
\end{equation}
where  $1 < \rho < |\Phi_E(\zeta)|.$
Combining \eqref{a} and \eqref{b}, we get
\begin{equation}\label{c}
\frac{1}{R_{\zeta}({\bf f},{\bf m})}=c\cdot  \limsup_{n\rightarrow\infty}\left|\int_{\Gamma_{\rho}}\frac{g(t)}{a_{n+1}(t)}\left(Q_{{\bf m}}^{\bf f}(t)-Q_{n,{\bf m}}(t)\right)dt\right|^{1/n}.
\end{equation}
This equality is impossible because from \eqref{eq5}, \eqref{eq20}, and \eqref{d} it is not hard to deduce that \eqref{c} is strictly less than ${1}/{R_{\zeta}({\bf f},{\bf m})}$.   This proves the equality in \eqref{eq1teo}.

\medskip

If $\xi$ is any one of the system poles of $\bf f$ and $\tau$ its order, from \eqref{eq14} and \eqref{eq17}, we have
\begin{equation}\label{eq21.2}
\max_{j=0\dots,l}\limsup_{n\rightarrow\infty}|Q_{n,{\bf m}}^{(j)}(\xi)|^{1/n}\leq\frac{|\Phi_E(\xi)|}{R_{\xi,l+1}({\bf f},{\bf m})},\quad l=0,1,\dots,\tau-1.
\end{equation}
Now we are ready to prove \eqref{eq6.1}. Let us fix $k\in\lbrace 1,\dots,d\rbrace$. Let $\mathcal{K}$ be a compact subset contained in $D_{k}^*({\bf f},{\bf m})\setminus\mathcal{P}_{\mathbf{m}}^{\mathbf{f}}$. Take $\delta>0$ sufficiently small so that
\[1<\rho:=R_{k}^*({\bf f},{\bf m})-\delta,\quad \quad  \mathcal{K} \subset D_\rho, \quad\quad   \bigcup_{j=1}^{N_k} \lbrace z\in \mathbb{C}:|z-\xi_j|\leq \delta\rbrace\subset D_{\rho}\setminus \mathcal{K}, \] where $\xi_1,\ldots,\xi_{N_k}$ are the poles of $f_k$ in $D_k^*({\bf f},{\bf m})$. Set
\[
C_j:=\{z\in \mathbb{C}:|z-\xi_j|=\delta\}.
\]
Let $\Gamma_{\rho,\delta}$ be the positively oriented curve determined by $\Gamma_\rho$ and those circles $C_j$. On account of Definition \ref{multipdef}, using Hermite's formula, we have
\begin{equation}\label{e}
(Q_{n,{\bf m}}f_k-P_{n,{\bf m},k})(z)=\frac{1}{2\pi i}\int_{\Gamma_{\rho,\delta}}\frac{a_{n+1}(z)}{a_{n+1}(t)}\frac{(Q_{n,{\bf m}}f_k)(t)}{t-z}dt.
\end{equation}
From \eqref{eq5} it readily follows that for all $z\in \mathcal{K},$
\begin{equation}\label{f}
\limsup_{n\rightarrow\infty} \left|\frac{1}{2\pi i}\int_{\Gamma_{\rho}}\frac{a_{n+1}(z)}{a_{n+1}(t)}\frac{(Q_{n,{\bf m}}f_k)(t)}{t-z}dt\right|^{1/n} \leq \frac{\|\Phi_E\|_{\mathcal{K}}}{R_{k}^*({\bf f},{\bf m})}.
\end{equation}
Let $\hat{\tau}_j$ be the order of $\xi_j$ as pole of $f_k$. Using the expansion
\[ Q_{n,{\bf m}}(t) = \sum_{l=0}^{|{\bf m}|} \frac{Q_{n,{\bf m}}^{(l)}(\xi_j)}{l!}(t - \xi_j)^{l},
\]
for the circle $C_j$ we have
\begin{equation}\label{g}
\frac{1}{2\pi i}\int_{C_j}\frac{a_{n+1}(z)}{a_{n+1}(t)}\frac{(Q_{n,{\bf m}}f_k)(t)}{t-z}dt = \sum\limits_{l=0}^{\hat{\tau}_j-1}\frac{1}{2\pi i}\int_{C_j} \frac{a_{n+1}(z)}{a_{n+1}(t)}\frac{(t - \xi_j)^{\hat{\tau}_j}f_k(t)}{(t-\xi_j)^{\hat{\tau}_j-l}}\frac{Q_{n,{\bf m}}^{(l)}(\xi_j)}{l!(t-z)}dt
\end{equation}
because the function under the integral sign is analytic inside $C_j$ for $\hat{\tau}_j \leq l \leq |{\bf m}|$. Now, \eqref{eq5} and \eqref{eq21.2} allow to deduce from \eqref{g} that for all $z\in \mathcal{K},$
\begin{equation}\label{h}
\limsup_{n\rightarrow\infty} \left|\frac{1}{2\pi i}\int_{C_j}\frac{a_{n+1}(z)}{a_{n+1}(t)}\frac{(Q_{n,{\bf m}}f_k)(t)}{t-z}dt\right|^{1/n} \leq \frac{\|\Phi_E\|_{\mathcal{K}}}{|\Phi_E(\xi_j)|}\frac{|\Phi_E(\xi_j)|}{R_{\xi_j,\hat{\tau}_j}({\bf f},{\bf m})}.
\end{equation}
Finally, \eqref{e}, \eqref{f}, and \eqref{h} give \eqref{eq6.1}. \hfill $\Box$

\medskip

A  slight variation of the arguments employed above allows to deduce the following corollary of independent interest.

\begin{Corollary} Let ${\bf f}\in \mathcal{H}(E)^d$ and fix a multi-index ${\bf m}\in \mathbb{N}^d.$ Suppose that \eqref{eq5} takes place and $\bf f$ has exactly $|\bf m|$ system poles with respect to $\bf m$. Then, for every system pole $\xi$ of $\bf f,$
\begin{equation}\label{eq21.1}
\max_{j=0\dots,l} \limsup_{n\rightarrow\infty}|Q_{n,{\bf m}}^{(j)}(\xi)|^{1/n}=\frac{|\Phi_E(\xi)|}{R_{\xi,l+1}({\bf f},{\bf m})},\quad l=0,1,\dots,\tau-1.
\end{equation}
where $\tau$ is the order of $\xi$.
\end{Corollary}

\begin{pf} If \eqref{eq21.1} fails, due to \eqref{eq21.2}, there is a system pole $\xi$ of $\bf f$ of order $\tau$ such that for some $l, 0 \leq l < \tau$
\begin{equation}\label{i} \max_{j=0\dots,l} \limsup_{n\rightarrow\infty}|Q_{n,{\bf m}}^{(j)}(\xi)|^{1/n}< \frac{|\Phi_E(\xi)|}{R_{\xi,l+1}({\bf f},{\bf m})}.
\end{equation}
Now, we argue by contradiction as in the proof of the equality in \eqref{eq1teo}.
\medskip

Choose a polynomial combination $g$ as in \eqref{eq21} that is analytic on a neighborhood of $\overline{D}_{|\Phi_E(\xi)|}$ except for a pole of order $s(\leq l+1)$ at $z=\xi$ with $\rho_s(g)= R_{\xi,l+1}({\bf f},{{\bf m}})$. Set $Q_{\bf m}^{\bf f}=Q_{{\bf m}}$. Take $\delta>0$ sufficiently small and $1<\rho < \rho_s(g)$. Let $\Gamma_{\rho,\delta}$ be the positively oriented curve determined by $\Gamma_\rho$ and $\lbrace t\in \mathbb{C}:|t-\xi|=\delta\rbrace$. Arguing as in \eqref{a}, it follows from \eqref{cauchyhadamard} that
\begin{equation}\label{eq21.3}
\frac{1}{\rho_s(g)}=c \cdot \limsup_{n\rightarrow\infty}\left|\int_{\Gamma_{\rho,\delta}}\frac{Q_{{\bf m}}(t)g(t)}{a_{n+1}(t)}dt\right|^{1/n}.
\end{equation}

\medskip

The function
\[
\frac{H_n(z)}{a_{n+1}(z)}=\frac{Q_{n,{\bf m}}(z)g(z)-\sum\limits_{k=1}^{d}p_{k}(z)P_{n,{\bf m},k}(z)}{a_{n+1}(z)}
\]
is analytic in $D_{\rho_s}(g)\setminus\lbrace\xi\rbrace$ and
\[
\int_{\Gamma_{\rho,\delta}}\frac{H_{n}(t)}{a_{n+1}(t)}dt=0.
\]
Set $P_n:=\sum\limits_{k=1}^{d}p_{k}P_{n,{\bf m},k}$ and $h:=(t-\xi)^sg$. Obviously,
\[
Q_{{\bf m}}g = (Q_{{\bf m}} -Q_{n,{\bf m}})g+P_n+H_n,
\]
and since $\deg P_n\leq n-1$, we obtain
\begin{align*}
&\int_{\Gamma_{\rho,\delta}}\frac{Q_{{\bf m}}(t)g(t)}{a_{n+1}(t)}dt=\int_{\Gamma_{\rho,\delta}}\frac{[Q_{{\bf m}}-Q_{n,{\bf m}}](t)h(t)}{(t-\xi)^sa_{n+1}(t)}dt\\
&=\int_{\Gamma_{\rho}}\frac{[Q_{{\bf m}}-Q_{n,{\bf m}}](t)h(t)}{(t-\xi)^sa_{n+1}(t)}dt-\sum_{j=0}^{|{\bf m}|}\int_{|t-\xi|=\delta}\frac{[Q_{{\bf m}}^{(j)}-Q_{n,{\bf m}}^{(j)}](\xi)h(t)}{j!(t-\xi)^{s-j}a_{n+1}(t)}dt\\
&=\int_{\Gamma_{\rho}}\frac{[Q_{{\bf m}}-Q_{n,{\bf m}}](t)h(t)}{(t-\xi)^sa_{n+1}(t)}dt+\sum_{j=0}^{s-1}\int_{|t-\xi|=\delta}\frac{Q_{n,{\bf m}}^{(j)}(\xi)h(t)}{j!(t-\xi)^{s-j}a_{n+1}(t)}dt.
\end{align*}
Estimating these integrals, using \eqref{eq5}, \eqref{eq1teo}, and the assumption \eqref{i}, it is easy to deduce that
\[ c \cdot \limsup_{n \to \infty}\left|\int_{\Gamma_{\rho,\delta}}\frac{Q_{{\bf m}}(t)g(t)}{a_{n+1}(t)}dt\right|^{1/n} < \frac{1}{\rho_s(g)}
\]
which contradicts \eqref{eq21.3}. Therefore, \eqref{i} cannot occur and there is equality in \eqref{eq21.1}. \end{pf}

\begin{remark} We wish to underline that for the proof of the previous results, excluding the equality in \eqref{eq1teo} and \eqref{eq21.1}, it would have been sufficient to assume that the table of points verifies \eqref{raiz} instead of \eqref{eq5}. The condition \eqref{eq5} has only been used in order to have the Cauchy Hadamard type formula \eqref{cauchyhadamard}. For the inverse type statement  $(b)\Rightarrow (a)$ the stronger assumption \eqref{eq5} is much more substantial.
\end{remark}

\section{Inverse statements}

\subsection{Some auxiliary results}

Let
\begin{equation}\label{eq22}
f(z)=\sum_{n=0}^{\infty}f_nz^n
\end{equation}
be a power series convergent in some neighborhood of the point $z=0$ whose radius of convergence we denote $R_0(f)$. We find it convenient to denote the $n$-th Taylor coefficient $f_n$ of the expansion $f$ also by $\left[f\right]_n$. According to the Cauchy-Hadamard formula
\[
R_{0}(f)=\left( \limsup_{n \rightarrow \infty} |[f]_n|^{1/n} \right)^{-1}.
\]
When $R_{0}(f) >0$ we define
\[
U_{\delta}(f):=\left\lbrace z \in \mathbb{C}: R_0(f)e^{-\delta}<|z|<R_0(f)e^{\delta}\right\rbrace.
\]

The following theorem was proved by V.I. Buslaev in \cite[Theorem 2]{Bus1}.

\medskip

\noindent
{\bf Buslaev's Theorem.} Suppose that $\delta>0$ and the power series \eqref{eq22} is such that  $0<R_0(f)<\infty$ and
\begin{equation}\label{eq23}
\left[f\alpha_n\right]_n=o(R_0(f)^{-n}e^{-n\delta}),
\end{equation}
where $\alpha_n\in \mathcal{H}(U_{\delta}(f))$ $(n=1,2,\dots)$ and $\lim_{n\rightarrow\infty}\alpha_n(z)=\alpha(z)$ $(z\in U_{\delta}(f))$. Then $\alpha$ has at least one zero on the circle $|z|=R_0(f)$, and the terms of the sequence $\lbrace f_n\rbrace$ $(n=1,2,\dots)$ satisfy
\begin{equation}\label{eq24}
f_{n+k}+ \beta_{n,1}f_{n+k-1} + \cdots + \beta_{n,k} f_{n} = 0,\qquad \lim_{n \rightarrow \infty} \beta_{n,p} = \beta_p, \qquad p=1,\ldots,k,
\end{equation}
the polynomial $\beta(z):= 1+\beta_1z + \cdots +\beta_k z^{k}$ $(\beta_k\neq 0)$ divides $\alpha$, all its roots are equal in modulus to $R_0(f)$, and at least one of them is a singular point of $f$.

\medskip

Buslaev's theorem can be supplemented by the following assertion (see \cite{Bus1}).

\medskip

\noindent
{\bf Supplement to Buslaev's Theorem.} Suppose that the power series \eqref{eq22} is not a polynomial, $R_0(f)=\infty$, and
\begin{equation}\label{eq25}
\alpha_{n,0}f_n+\alpha_{n,-1}f_{n+1}+\cdots=0\quad (n=1,2,\dots)
\end{equation}
where the $\alpha_n(z)=\sum_{p=0}^{\infty}\alpha_{n,-p}z^{-p}$ $(n=1,2,\dots)$ are holomorphic and converge to $\alpha(z)$ in the exterior of some disk as $n\rightarrow\infty$. Then $\alpha(\infty)=0$, and the coefficients $\lbrace f_n\rbrace$ of the series \eqref{eq22} satisfy
\begin{equation}\label{eq26}
\epsilon_{n,0}f_{n} + \cdots + \epsilon_{n,-N+1}f_{n+N-1}+f_{n+N} = 0,\qquad \lim_{n \rightarrow \infty} \epsilon_{n,p} = \epsilon_p, \qquad p=0,-1,\ldots,-N+1,
\end{equation}
$N$ being the multiplicity of the zero of $\alpha$ at $z=\infty$.

\medskip

This result will be useful in the next section to prove Lemma \ref{lemma1}.

\subsection{Incomplete multipoint Pad\'{e} approximants}

Let us introduce the notion of incomplete multipoint Pad\'{e} approximants.  A similar concept turned out to be effective in the study of Hermite-Pad\'{e} approximation in \cite{cacoq1} and \cite{cacoq2} for proving results of inverse type.
\begin{Definition}\label{multipincomp}\textup{
Let $f\in \mathcal{H}(E)$. Fix $m\geq m^*\geq 1$ and $n\geq m$. We say that the rational function $R_{n,m}$ is an \emph{incomplete multipoint Pad\'{e} approximant of type $(n,m,m^*)$ corresponding to $f$} if $R_{n,m}$ is the quotient of any two polynomials $P_{n,m}$, $Q_{n,m}$ that verify
\begin{enumerate}
\item[c.1)] $\deg P_{n,m}\leq n-m^*$, $\deg Q_{n,m}\leq m$, $Q_{n,m}\not\equiv 0$,
\item[c.2)] $\displaystyle \frac{Q_{n,m}f-P_{n,m}}{a_{n+1}}\in \mathcal{H}(E),$
\end{enumerate}
where $a_n(z)=\prod_{k=1}^{n}(z-\alpha_{n,k})$.}
\end{Definition}

Since $Q_{n,m}\not\equiv 0$, we normalize it to be monic. We call $Q_{n,m}$ the denominator of the corresponding $(n,m,m^*)$ incomplete multipoint Pad\'{e} approximant of $f$. Notice that for each $k=1,\dots,d$, the polynomial $Q_{n,{\bf m}},$ given in Definition \ref{multipdef}, is a denominator of an $(n,|{\bf m}|,m_k)$ incomplete multipoint Pad\'{e} approximant of $f_k$.

\medskip

In this section, we will study the relation between the convergence of $Q_{n,m}$ and some analytic properties of $f$.

\begin{Lemma}\label{lemma1}
Let $f\in\mathcal{H}(E)$ and fix $m\geq m^*\geq 1$. Suppose that $f$ is not a rational function with at most $m^*-1$ poles and there exists a polynomial $Q_m$ of degree $m$ such that
\begin{equation}\label{eqlemma1}
\limsup_{n\rightarrow\infty}\|Q_{n,m}-Q_m\|^{1/n}\leq\theta<1.
\end{equation}
Then, either $f$ has exactly $m^*$ poles in $D_{\rho_{m^*}(f)}$ or $\rho_0(Q_mf)>\rho_{m^*}(f)$, where $\rho_{m^*}(f)$ is the index of the largest canonical region to which $f$ can be extended as a meromorphic function with at most $m^*$ poles counting multiplicities.
\end{Lemma}

\noindent {\bf Proof.}
Let $\lbrace\xi_1,\dots,\xi_{\omega}\rbrace$ be the distinct poles of $f$ in $D_{\rho_{m^*}(f)}$ and $\tau_1,\dots,\tau_{\omega}$ be their orders, respectively. Consequently,
\[
\sum_{j=1}^{\omega}\tau_j\leq m^*.
\]
Modifying conveniently the proof of \eqref{eq induc}, one can show that for $j=1,\dots,\omega$
\begin{equation}\label{eq27}
\limsup_{n\rightarrow\infty}|Q_{n,m}^{(\nu)}(\xi_j)|^{1/n}\leq\frac{|\Phi_E(\xi_j)|}{\rho_{m^*}(f)}<1, \qquad \nu=0,1,\dots,\tau_j-1.
\end{equation}
Since the sequence of polynomials $Q_{n,m}$ converges to $Q_m$, \eqref{eq27} entails that $\xi_j$ is a zero of $Q_m$ of multiplicity at least $\tau_j$. Being this the case,  we have
\[
\rho_0(Q_mf)\geq\rho_{m^*}(f).
\]

\medskip

Suppose that $\rho_0(Q_mf)=\rho_{m^*}(f)$. To conclude the proof, let us show that in this situation $f$ has exactly $m^*$ poles in $D_{\rho_{m^*}(f)}$.
To the contrary, suppose that $f$ has in $D_{\rho_{m^*}(f)}$ at most $m^*-1$ poles. Then, there exists a polynomial $\deg Q_{m^*}<m^*$ such that
\[
\rho_0(Q_{m^*}f)=\rho_{m^*}(f)=\rho_0(Q_mQ_{m^*}f).
\]
It follows from Definition \ref{multipincomp} that
\[
\frac{Q_{m^*}(Q_{n,m}f-P_{n,m})}{a_{n+1}}\in \mathcal{H}(E).
\]
Then
\begin{equation*}
\int_{\Gamma_{\rho}}\frac{Q_{m^*}(z)(Q_{n,m}f-P_{n,m})(z)}{a_{n+1}(z)}dz=0,
\end{equation*}
where $1<\rho<\rho_{m^*}(f)$. Since each one of the $n+1$ zeros of the polynomial $a_{n+1}$ lies on $E$ and $\deg (Q_{m^*}P_{n,m})\leq n-1$, it follows that
\begin{equation*}
\int_{\Gamma_{\rho}}\frac{Q_{m^*}(z)P_{n,m}(z)}{a_{n+1}(z)}dz=0.
\end{equation*}
Therefore,
\begin{equation}\label{eq28}
\int_{\Gamma_{\rho}}\frac{Q_{m^*}(z)Q_{n,m}(z)f(z)}{a_{n+1}(z)}dz=0.
\end{equation}
 Then, by \eqref{cauchyhadamard},
\begin{align*}
\frac{1}{\rho_{m^*}(f)}&=\frac{1}{\rho_0(Q_mQ_{m^*}f)}=c\cdot \limsup_{n\rightarrow\infty}\left|\int_{\Gamma_{\rho}}\frac{(Q_mQ_{m^*}f)(t)}{a_{n+1}(t)}dt\right|^{1/n}\\
&=c\cdot \limsup_{n\rightarrow\infty}\left|\int_{\Gamma_{\rho}}\frac{(Q_{m^*}f)(t)}{a_{n+1}(t)}\left(Q_{n,m}-Q_m\right)(t)dt\right|^{1/n}.
\end{align*}
Using \eqref{eq5}  and \eqref{eqlemma1} to estimate the last integral, it readily follows that
 \[
\frac{1}{\rho_{m^*}(f)}\leq\frac{\theta}{\rho_{m^*}(f)},  \qquad\theta<1,
\]
which implies that $\rho_{m^*}(f)=\infty$. Now, let us show that this is not possible.

\medskip

Take $F(w):=Q_{m^*}(\Psi_E(w))f(\Psi_E(w))$, where $\Psi_E=\Phi_E^{-1}$. Let $\gamma$ be a contour encircling $\{w\in \mathbb{C}:|w|=1\}$ lying in the domain of holomorphy of $F$. Using \eqref{eq28}, we obtain
\[
0 = \int_{\gamma}\frac{F(w)Q_{n,m}(\Psi_E(w)) }{a_{n+1}(\Psi_E(w))}\Psi_E'(w)dw =
\]
\[  \int_{\gamma} F(w) \frac{Q_{n,m}(\Psi_E(w)) }{w^{m}}\frac{w^{n+1}}{a_{n+1}(\Psi_E(w))}\Psi_E'(w)\frac{dw}{w^{n+1-m}}
\]
Setting
\[ \alpha_n(w)=\frac{Q_{n,m}(\Psi_E(w)) }{w^{m}}\frac{(cw)^{n+1}}{a_{n+1}(\Psi_E(w))}\Psi_E'(w),
\]
the previous equality means that
\begin{equation}\label{eq29}
\left[F\alpha_n\right]_{n-m}=0.
\end{equation}
The functions $\alpha_n$ $(n=1,2,\dots)$ are holomorphic in the exterior of the unit disk (including $w=\infty$) and, due to \eqref{eq5} and \eqref{eqlemma1}, converge as $n\rightarrow\infty$ to
\[
\alpha(w)=\Psi_E'(w)\frac{ Q_m(\Psi_E(w))}{w^{m}G(\Psi_E(w))}=\sum_{p=0}^{\infty}\alpha_{-p}w^{-p}, \qquad \alpha_0 = \alpha(\infty) \neq 0.
\]

Let $\sum_{n=-\infty}^{\infty}F_nw^{n}$ be the Laurent expansion of the function $F$ outside the unit circle, i.e:
\[
F(w) =\sum_{n=-\infty}^{\infty}F_nw^{n}=F_1(w)+F_2(w),
\]
where $F_1(w)=\sum_{n=0}^{\infty}F_nw^{n}$. Then, $R_0(F_1) = \infty$ and \eqref{eq29} holds (for all sufficiently large $n$) replacing $F$ with $F_1$. According to the Supplement to Buslaev's Theorem and the fact that $\alpha(\infty)\neq 0$, we get that $F_1$ must be a polynomial. Consequently, $F$ is either analytic or has a pole at $\infty$. In turn this implies that $Q_{m^*}f$ is either analytic or has a pole at $\infty$. However, $Q_{m^*}f$ is an entire function because it is holomorphic in $\mathbb{C}$ since $R_0(Q_{m^*}f) = \infty$. Therefore, $Q_{m^*}f$ is a polynomial, or what is the same $f$ is a rational function with at most $m^* -1$ poles against our hypothesis on $f$. This contradiction implies that the assumption that  $f$ had in  $D_{\rho_{m^*}(f)}$ at most $m^* -1$ poles is impossible. So the number of poles on $f$ in $D_{\rho_{m^*}(f)}$ must equal $m^*$. \hfill $\Box$

\subsection{Polynomial independence}

Let us introduce the concept of polynomial independence of a vector of functions.

\begin{Definition}\label{indpoly}\textup{
A vector $\mathbf{f}=(f_1,\dots,f_d)\in \mathcal{H}(E)^{d}$ is said to be \emph{polynomially independent with respect to $\mathbf{m}=(m_1,\dots,m_d)\in\mathbb{N}^{d}$} if there do not exist polynomials $p_1,\dots,p_d$, at least one of which is non-null, such that
\begin{enumerate}
\item[(i)] $\deg p_k< m_k,$ $k=1,\dots,d$,
\item[(ii)] $\sum_{k=1}^dp_kf_k$ is a polynomial.
\end{enumerate}}
\end{Definition}

In particular, polynomial independence implies that for each $k=1,\dots,d$, $f_k$ is not a rational function with at most $m_k-1$ poles.

\begin{Lemma}\label{lemma2}
Let $\mathbf{f}\in\mathcal{H}(E)^d$ and fix a multi-index $\mathbf{m}\in\mathbb{N}^{d}$. Suppose that for all $n\geq n_0$, the polynomial $Q_{n,\mathbf{m}}$ is unique and $\deg Q_{n,\mathbf{m}}=|\mathbf{m}|$. Then the system $\mathbf{f}$ is polynomially independent with respect to $\mathbf{m}$.
\end{Lemma}

\noindent {\bf Proof.} Except for a small detail, the proof coincides with that of \cite[Lemma 3.2]{cacoq2}. Given $\mathbf{f}:=(f_1,\dots,f_d)\in \mathcal{H}(E)^d$ and $\mathbf{m}:=(m_1,\dots,m_d)\in\mathbb{N}^{d}$, we consider the associated system
\begin{equation*}
\mathbf{\overline{f}}:=(f_1,\dots,z^{m_1-1}f_1,f_2,\dots,z^{m_d-1}f_d)=(\overline{f}_1,\dots,\overline{f}_{{|\bf m}|}).
\end{equation*}
We also define an associated multi-index $\overline{\bf m}:=(1,\dots,1)$ with $|{\bf m}|=|\overline{\bf m}|$. The systems ${\bf f}$ and $\mathbf{\overline{f}}$ share most properties. In particular, poles and system poles  of $({\bf f}.{\bf m})$ and $(\mathbf{\overline{f},\overline{\bf m}})$ coincide and $\mathbf{f}$ is polynomially independent with respect to $\mathbf{m}$ if and only if $\overline{\mathbf{f}}$ is polynomially independent with respect to $\overline{\mathbf{m}}$. Passing to $(\mathbf{\overline{f},\overline{\bf m}})$ if necessary and relabeling the functions,   we can assume without loss of generality that ${\bf m}=(1,\dots,1)$ and $d=|{\bf m}|$.

\medskip

Suppose that there exist constants $c_k$, $k=1,\dots,d$, not all zero, such that $\sum_{k=1}^dc_kf_k$ is a polynomial. Without loss of generality, we can assume that $c_1\neq 0$. Then,
\[
f_1=p-\sum_{k=2}^dc_kf_k,
\]
where $p$ is a polynomial of degree $N$.

\medskip

On the other hand, for each $n\geq d-1$, there exist polynomials $Q_{n}$, $P_{n,k}$, $k=2,\dots,d$, such that for all $k=2,\dots,d$,

\begin{enumerate}
\item[-] $\deg P_{n,k}\leq n-1$, $\deg Q_{n}\leq d-1$, $Q_{n}\not\equiv 0$,
\item[-] $ \displaystyle \frac{Q_{n}f_k-P_{n,k}}{a_{n+1}}\in \mathcal{H}(E)$.
\end{enumerate}
Therefore,
\[
\frac{Q_{n}\left( p-\sum_{k=2}^dc_kf_k\right)-\left( Q_{n}p-\sum_{k=2}^dc_kP_{n,k}\right)}{a_{n+1}}\in \mathcal{H}(E)
\]
and, for $n\geq d+N$, the polynomial $P_{n,1}=Q_{n}p-\sum_{k=2}^dc_kP_{n,k}$ verifies $\deg P_{n,1}\leq n-1$. Thus, for all $n$ sufficiently large, the polynomials $P_{n,k}$, $k=1,\dots,d$ satisfy Definition \ref{multipdef} with respect to ${\bf f}$ and ${\bf m}$. Naturally, $Q_n$ gives rise to a polynomial $Q_{n,{\bf m}}$ with $\deg Q_{n,{\bf m}}< d=|{\bf m}|$ against our assumption on $Q_{n,{\bf m}}$. \hfill $\Box$

\medskip

The following corollary is a straightforward consequence of Lemma \ref{lemma1}.
\begin{Corollary}\label{corollary1}
Let $\mathbf{f}\in\mathcal{H}(E)^d$ and fix a multi-index $\mathbf{m}\in\mathbb{N}^{d}$. Assume that $\mathbf{f}$ is polynomially independent with respect to $\mathbf{m}$ and there exists a polynomial $Q_{\mathbf{m}}$ of degree $|\mathbf{m}|$ such that.
\begin{equation}\label{eqcorollary1}
\limsup_{n\rightarrow\infty}\|Q_{n,\mathbf{m}}-Q_{\mathbf{m}}\|^{1/n}\leq\theta<1.
\end{equation}
Then for each $k=1,\dots,d$, either $f_k$ has exactly $m_k$ poles in $D_{\rho_{m_k}(f_k)}$ or $\rho_0(Q_{\mathbf{m}}f_k)>\rho_{m_k}(f_k)$.
\end{Corollary}

Before proving the inverse statement of Theorem \ref{maintheorem}, we wish to describe some properties of system poles. For the proof see \cite[Lemma 3.5]{cacoq2}.

\begin{Lemma}\label{lemma7}
Let $\mathbf{f}\in\mathcal{H}(E)^d$ and $\mathbf{m}\in\mathbb{N}^{d}$. Then, $\mathbf{f}$ can have at most $|\mathbf{m}|$ system poles with respect to $\mathbf{m}$ (counting their order). Moreover, if the system $\mathbf{f}$ has exactly $|\mathbf{m}|$ system poles with respect to $\mathbf{m}$ and $\xi$ is a system pole of order $\tau$, then for all $s>\tau$ there can be no polynomial combination of the form \eqref{eqsystempole2} holomorphic in a neighborhood of $\overline{D}_{|\Phi(\xi)|}$ except for a pole at $z=\xi$ of exact order $s$.
\end{Lemma}
\subsection{Proof $(b)\Rightarrow(a)$}

Arguing as in the beginning of the proof of Lemma \ref{lemma2}, we can assume without loss of generality that ${\bf m} = (1,1,\ldots,1)$ and $|{\bf m}| = d$. From Definition \ref{indpoly}, it readily follows that ${\bf f}$ is polynomially independent with respect to ${\bf m}$ if and only if there do not exist constants $c_k$, $k=1,\dots,|{\bf m}|$, not all zero, such that
\[
\sum_{k=1}^{|{\bf m}|}c_k {f}_k
\]
is a polynomial. Due to Lemma \ref{lemma2}, on account of the hypothesis, we know that $\mathbf{f}$ is polynomially independent with respect to $\mathbf{m}$. We must show that ${\bf f}$ has exactly $|{ {\bf m}}|$ system poles with respect to ${{\bf m}}.$

\medskip

The auxiliary results that we have obtained allow us to adapt the proof employed in \cite{cacoq2}, where classical Hermite-Pad\'e approximation was considered.  For completeness we include the whole proof.

\medskip

The scheme is as follows. First, we collect a set of $|{\bf  {m}}|$ candidates to be system poles of ${\mathbf{ {f}}}$ (counting their orders) and prove that they are zeros of $Q_{\mathbf{m}}$. In the second part we prove that all these points previously selected are actually system poles of ${\mathbf{ {f}}}$.

\medskip

Notice that for each $k=1,\dots,|{\bf  {m}}|$, by Corollary \ref{corollary1}, either $D_{\rho_1( {f}_k)}$ contains exactly one pole of $ {f}_k$ and it is a zero of $Q_{\mathbf{m}}$, or $\rho_0(Q_{\mathbf{m}} {f}_k)>\rho_1( {f}_k)$. Hence, $D_{\rho_0({\mathbf{ {f}}})}\neq\mathbb{C}$ and $Q_{\mathbf{m}}$ contains as zeros all the poles of $ {f}_k$ on the boundary of $D_{\rho_0( {f}_k)}$ counting their order for $k=1,\dots,|\mathbf{ {m}}|$. Moreover, the function $ {f}_k$ cannot have on the boundary of $D_{\rho_0( {f}_k)}$ singularities other than poles. Hence, the poles of ${\mathbf{ {f}}}$ on the boundary of $D_{\rho_0({\mathbf{ {f}}})}$ are all zeros of $Q_{\mathbf{m}}$ counting multiplicities and the boundary contains no other singularity except poles. Let us call them candidate system poles of ${\mathbf{ {f}}}$ and denote them by $a_1,\dots,a_{n_1}$ repeated according to their order. They constitute the first layer of candidate system poles of ${\mathbf{ {f}}}$.

\medskip

Since $\deg Q_{\mathbf{m}}=|\mathbf{ {m}}|$, $n_1\leq |\mathbf{ {m}}|$. If $n_1=|\mathbf{ {m}}|$, we are done. Let us assume that $n_1<|\mathbf{ {m}}|$ and let us find coefficients $c_1,\dots,c_{|\mathbf{ {m}}|}$ such that
\[
\sum\limits_{k=1}^{|\mathbf{ {m}}|}c_k {f}_k
\]
is holomorphic in a neighborhood of $\overline{D}_{\rho_0({\mathbf{ {f}}})}$. Finding those $c_1,\dots,c_{|\mathbf{ {m}}|}$ reduces to solving a homogeneous system of $n_1$ linear equations with $|\mathbf{ {m}}|$ unknowns. In fact, if $z=a$ is a candidate system pole of ${\mathbf{ {f}}}$ with order $\tau$, we obtain $\tau$ equations choosing the coefficients $c_k$ so that
\begin{equation}\label{eq30}
\int_{|\omega-a|=\delta}(\omega-a)^k\left(\sum\limits_{k=1}^{|\mathbf{ {m}}|}c_k {f}_k(\omega)\right)d\omega=0,\qquad k=0,\dots,\tau-1.
\end{equation}
We obtain the same type of   equations for each distinct candidate system pole on the boundary of $D_{\rho_0({\mathbf{ {f}}})}$. Combining these equations, we obtain a homogeneous system of $n_1$ linear equations with $|\mathbf{ {m}}|$ unknowns. Moreover, this homogeneous system of linear equations has at least $|\mathbf{ {m}}|-n_1$ linearly independent solutions, which we denote by $\mathbf{c}_j^1$, $j=1,\dots,|\mathbf{m}|-n_1^*$, where $n_1^*\leq n_1$ denotes the rank of the system of equations.

\medskip

Let
\[
\mathbf{c}_j^1:=(c_{j,1}^1,\dots,c_{j,|\mathbf{ {m}}|}^1),\qquad j=1,\dots,|\mathbf{ {m}}|-n_1^*.
\]
Define the $(|\mathbf{ {m}}|-n_1^*)\times |\mathbf{ {m}}|$ dimensional matrix
\[
 C^1:=\begin{pmatrix} \mathbf{c}_1^1 \\ \vdots \\ \mathbf{c}_{|\mathbf{m}|-n_1^*}^1 \end{pmatrix}.
\]
Define the vector $\mathbf{g}_1$ of $|\mathbf{m}|-n_1^*$ functions given by
\[
\mathbf{g}_1^t:=C^1{\mathbf{ {f}}}^t=(g_{1,1},\dots,g_{1,|\mathbf{m}|-n_1^*})^t,
\]
where $(\cdot)^t$ means taking transpose. Since all the rows of $C^1$ are non-null and ${\mathbf{ {f}}}$ is polynomially independent with respect to ${\mathbf{ {m}}}$, none of the functions
\[
g_{1,j}=\sum\limits_{k=1}^{|\mathbf{ {m}}|}c_{j,k}^1 {f}_k,\qquad j=1,\dots,|\mathbf{ {m}}|-n_1^*,
\]
are polynomials.

\medskip

Consider the canonical domain
\[
D_{\rho_0(\mathbf{g}_1)}=\bigcap\limits_{j=1}^{|\mathbf{m}|-n_1^*}D_{\rho_0(g_{1,j})}.
\]
Obviously, by construction, $D_{\rho_0({\mathbf{ {f}}})}$ is strictly included in $D_{\rho_0(\mathbf{g}_1)}$. Therefore, for each $j=1,\dots,|\mathbf{ {m}}|-n_1^*$, $Q_{n,\mathbf{m}}$ is a denominator of an $(n,|\mathbf{ {m}}|,1)$ multipoint incomplete Pad\'{e} approximant of $g_{1,j}$. Since the $g_{1,j}$ are not polynomials, by Lemma \ref{lemma1} with $m^*=1$, for each $j=1,\dots,|\mathbf{ {m}}|-n_1^*$, either $D_{\rho_1(g_{1,j})}$ contains exactly one pole of $g_{1,j}$ and it is a zero of $Q_{\mathbf{m}}$, or $\rho_0(Q_{\mathbf{m}}g_{1,j})>\rho_{1}(g_{1,j})$. In particular, $D_{\rho_0(\mathbf{g_1})}\neq\mathbb{C}$ and all the singularities of $\mathbf{g_1}$ on the boundary of $D_{\rho_0(\mathbf{g_1})}$ are poles which are zeros of $Q_{\mathbf{m}}$ counting their order. They form the next layer of candidate system poles of ${\mathbf{ {f}}}$.

\medskip

Denote by $a_{n_1+1},\dots,a_{n_1+n_2}$ the new candidate system poles. We repeat the arguments employed above. If $n_1+n_2=|\mathbf{ {m}}|$, we are done. Otherwise, $n_2<|\mathbf{ {m}}|-n_1\leq |\mathbf{ {m}}|-n_1^*$ and we eliminate the $n_2$ poles $a_{n_1+1},\dots,a_{n_1+n_2}$ as we did on the first layer. We have $|\mathbf{ {m}}|-n_1^*$ functions which are holomorphic on $D_{\rho_0(\mathbf{g_1})}$ and meromorphic on a neighborhood of $\overline{D}_{\rho_0(\mathbf{g_1})}$. The corresponding homogeneous system of linear equations, similar to \eqref{eq30}, has at least $|\mathbf{ {m}}|-n_1^*-n_2^*$ linearly independent solutions ${\bf c}_j^2$, where $n_2^*\leq n_2$ is the rank of the new system. Let
\[
\mathbf{c}_j^2:=(c_{j,1}^2,\dots,c_{j,|\mathbf{ {m}}|-n_1^*}^2),\qquad j=1,\dots,|\mathbf{ {m}}|-n_1^*-n_2^*.
\]
Define the $(|\mathbf{ {m}}|-n_1^*-n_2^*)\times (|\mathbf{ {m}}|-n_1^*)$ dimensional matrix
\[
 C^2:=\begin{pmatrix} \textbf{c}_1^2 \\ \vdots \\ \textbf{c}_{|\mathbf{m}|-n_1^*-n_2^*}^2 \end{pmatrix}.
\]
Define the vector $\mathbf{g}_2$ of $|\mathbf{ {m}}|-n_1^*-n_2^*$ functions given by
\[
\mathbf{g}_2^t:=C^2\mathbf{g}_1^t=C^2C^1{\mathbf{ {f}}}^t=(g_{2,1},\dots,g_{2,|\mathbf{ {m}}|-n_1^*-n_2^*})^t.
\]
It is a basic fact from linear algebra that if $C^1$ has full rank and $C^2$ has full rank, then $C^2C^1$ has full rank. This means that the rows of $C^2C^1$ are linearly independent, particulary, they are non-null. Therefore, none of the component functions of $\mathbf{g}_2$ are polynomials because of the polynomial independence of ${\mathbf{ {f}}}$ with respect to ${\mathbf{ {m}}}$. Thus, we can apply again Lemma \ref{lemma1}. Using finite induction, we find a total on $|{\mathbf{ {m}}}|$ candidate system poles.

\medskip

In fact, on each layer of system poles, $n_k\geq 1$. Therefore, in a finite number of steps, say $N-1$, their sum equals to $|{\mathbf{ {m}}}|$. Consequently, the number of candidate system poles of ${\mathbf{ {f}}}$ in some canonical domain, counting multiplicities, is exactly equal to $|{\mathbf{ {m}}}|$, and they are precisely the zeros of $Q_{\mathbf{m}}$ as we wanted to prove.

\medskip

Summarizing, in the $N-1$ steps we have taken, we have produced $N$ layers of candidate system poles. Each layer contains $n_k$ candidates, $k=1,\dots, N$. At the same time, on each step $k$, $k=1,\dots, N-1$, we have solved a system of $n_k$ linear equations, of rank $n_k^*$, with $|\mathbf{ {m}}|-n_1^*-\dots-n_k^*$, $n_k^*\leq n_k$, linearly independent solutions. We find ourselves on the $N$-th layer with $n_N$ candidates.

\medskip

Let us try to eliminate the poles on the last layer. Write the corresponding homogeneous system of linear equations as in \eqref{eq30}, and we get $n_N$ equations where
\[
n_N=|\mathbf{ {m}}|-n_1-\dots-n_{N-1}\leq |\mathbf{ {m}}|-n_1^*-\dots-n_{N-1}^*=:\overline{n}_N
\]
with $\overline{n}_N$ unknowns. For each candidate system pole $a$ of multiplicity $\tau$ on the $N$-th layer, we impose the equations
\begin{equation}\label{eq31}
\int_{|\omega-a|=\delta}(\omega-a)^j\left(\sum\limits_{k=1}^{\overline{n}_N}c_kg_{N-1,k}(\omega)\right)d\omega=0,\quad j=0,\dots,\tau-1,
\end{equation}
where $\delta$ is sufficiently small and the $g_{N-1,k}$, $k=1,\dots,\overline{n}_N$, are the functions associated with the linearly independent solutions produced on step $N-1$.

\medskip

Let $n_N^*$ be the rank of this last homogeneous system of linear equations. Assume that $n_k^*<n_k$ for some $k\in\lbrace 1,\dots, N\rbrace$. Then, the rank of the last system of equations is strictly less than the number of unknowns, namely $n_N^*<\overline{n}_N$. Therefore, repeating the same process, there exists a vector of functions
\[
\mathbf{g}_N:=(g_{N,1},\dots,g_{N,|\mathbf{ {m}}|-n_1^*-\dots-n_N^*})
\]
such that none of the $g_{N,k}$ is a polynomial because of the polynomial independence of ${\mathbf{{f}}}$ with respect to ${\mathbf{{m}}}$. Applying Lemma \ref{lemma1}, each $g_{N,k}$ has on the boundary of its canonical domain of analyticity a pole which is a zero of $Q_{\mathbf{m}}$. However, this is impossible because all the zeros of $Q_{\mathbf{m}}$ are strictly contained in that canonical domain. Consequently, $n_k^*=n_k$, for all $k= 1,\dots, N$.

\medskip

We conclude that all the $N$ homogeneous systems of linear equations that we have solved have full rank. This implies that if in any one of those $N$ systems of equations we equate one of its equations to $1$ instead of zero (see \eqref{eq30} or \eqref{eq31}), the corresponding nonhomogeneous system of linear equations has a solution. By the definition of a system pole, this implies that each candidate system pole is indeed a system pole of order at least equal to its multiplicity as zero of $Q_{\mathbf{m}}$. Moreover, by Lemma \ref{lemma7}, ${\mathbf{ {f}}}$ can have at most $|\mathbf{ {m}}|$ system poles with respect to ${\mathbf{m}}$; therefore, all candidate system poles are system poles, and their order coincides with the multiplicity of that point as a zero of $Q_{\mathbf{m}}$. This also means that $Q_{\mathbf{m}}= Q_{\mathbf{m}}^{\mathbf{f}}$. Thus, the proof  of the inverse type result is complete. \hfill $\Box$

{\noindent Nattapong Bosuwan\\
Department of Mathematics, Faculty of Science, Mahidol University,\\
Rama VI Road, Ratchathewi District,Bangkok 10400, Thailand\\
and\\
Centre of Excellence in Mathematics, CHE,\\
Si Ayutthaya Road, Bangkok 10400, Thailand\\}
{Email: nattapong.bos@mahidol.ac.th}\\

{\noindent Guillermo L\'{o}pez Lagomasino and Yanely Zaldivar Gerpe \\
Departament of Mathematics, Universidad Carlos III de Madrid, \\
Avda. Universidad 30, 28911 Legan\'{e}s, Madrid, Spain.\\}
{Emails: lago@math.uc3m.es. and yzaldiva@math.uc3m.es.}


\begin{thebibliography}{99}

\bibitem{Bos1} N. Bosuwan and G. L\'opez Lagomasino. Determining system poles using row sequences of orthogonal Hermite-Pad\'{e} approximants. J. Approx. Theory \textbf{231} (2018), 15-40.

\bibitem{Bos2} N. Bosuwan and G. L\'opez Lagomasino. Direct and inverse results on row sequences of simultaneous Padé-Faber approximants. Accepted in Medit. J. of Math.. arxiv 1801.03004

\bibitem{Bus1} {V.I. Buslaev}. Relations for the coefficients, and singular points of a function. Math. USSR Sb. {\bf 59} (1988), 349-377.

\bibitem{cacoq1} {J. Cacoq, B. de la Calle Ysern, and G. L\'{o}pez Lagomasino}. Incomplete Pad\'{e} approximation and convergence of row sequences of Hermite-Pad\'{e} approximants. J. Approx. Theory {\bf 170} (2013), 59-77.

\bibitem{cacoq2} {J. Cacoq, B. de la Calle Ysern, and G. L\'{o}pez Lagomasino}. Direct and inverse results on row sequences of Hermite-Pad\'{e} approximants. Constr. Approx. {\bf 38} (2013), 133-160.


\bibitem{gon1} {A.A. Gonchar}. On convergence of Pad\'e approximants for some classes of meromorphic  functions. Math. USSR Sb.  \textbf{26} (1975), 555-575.

\bibitem{gon3} A.A. Gonchar. Rational approximation of analytic functions. Proc. Steklov Inst. Math. {\bf 272} (2011), S44-S57.


\bibitem{gon2} A.A. Gonchar.  Poles of rows of the Pad\'e table and meromorphic continuation of functions. Sb. Math. \textbf{43} (1982), 527-546.

\bibitem{graves} {P.R. Graves-Morris and E.B. Saff}. A de Montessus theorem for vector-valued rational interpolants. Lecture Notes in Math., Vol. 1105, pp. 227-242, Springer, Berlin, (1984).

\bibitem{Mon} R. de Montessus de Ballore.  Sur les fractions continues alg\'ebriques. Bull. Soc. Math. France \textbf{30} (1902), 28-36.

\bibitem{Sidi} A. Sidi. A de Montessus type convergence study of a least-squares vector-valued rational interpolation procedure II. Comput. Methods Funct. Theory {\bf 10} (2010), 223-247.

\bibitem{VanBarel}  M. Van Barel and A. Bultheel. A new approach to the rational interpolation problem: the vector case. J. Comput.
Appl. Math. \textbf{33} (1990) 331-346.

\bibitem{Walsh} {J.L. Walsh}. Interpolation and Approximation by Rational Functions in the Complex Domain. 5th Ed.
Colloq. Publ. Vol. XX, Amer. Math. Soc., Providence, R. I. (1969).























\end{thebibliography}
\end{document}